\documentclass[12pt,leqno]{article}
\usepackage{amssymb}
\usepackage[mathscr]{eucal}
\usepackage{amsmath,amssymb,latexsym,theorem,bbm} %
\usepackage{color,url}
\usepackage{appendix}

\newcommand{\NN}{\mathbb{N}}

\newcommand{\RR}{\mathbb{R}}

\newcommand{\ZZ}{\mathbb{Z}}

\newcommand{\bM}{{\boldsymbol{M}}}

\newcommand{\bY}{{\boldsymbol{Y}}}

\newcommand{\cD}{{\mathcal D}}

\newcommand{\cF}{{\mathcal F}}

\newcommand{\EE}{\operatorname{\mathbb{E}}}
\newcommand{\PP}{\operatorname{\mathbb{P}}}
\newcommand{\OO}{\operatorname{O}}
\newcommand{\oo}{\operatorname{o}}

\newcommand{\vare}{\varepsilon}

\renewcommand{\mid}{\,|\,}

\renewcommand{\leq}{\leqslant}
\renewcommand{\geq}{\geqslant}

\newcommand{\distre}{\stackrel{\cD}{=}}

\newcommand{\proofend}{\hfill\mbox{$\Box$}}

\numberwithin{equation}{section}

\theoremstyle{change} \theorembodyfont{\em}
\newtheorem{Thm}{Theorem.}[section]
\newtheorem{Lem}[Thm]{Lemma.}
\newtheorem{Pro}[Thm]{Proposition.}

\newtheorem{Def}[Thm]{Definition.}

\theorembodyfont{\rm}
\newtheorem{Rem}[Thm]{Remark.}

\begin{document}

\begin{center}
 {\bfseries\Large {Regularly varying nonstationary second-order \\[2mm]
                   Galton--Watson processes with immigration}} \\[6mm]
 {\sc\large Zsuzsanna $\text{B\H{o}sze}^{*}$,
            \ Gyula $\text{Pap}^{*,\diamond}$}
\end{center}

\vskip0.2cm

\noindent
 * Bolyai Institute, University of Szeged,
   Aradi v\'ertan\'uk tere 1, H-6720 Szeged, Hungary.

\noindent e-mail: Bosze.Zsuzsanna@stud.u-szeged.hu (Zs. B\H{o}sze),
                  papgy@math.u-szeged.hu (G. Pap).

\noindent $\diamond$ Corresponding author.



\renewcommand{\thefootnote}{}
\vspace*{0.2cm}
\footnote{Supported by the Hungarian Croatian Intergovernmental S \& T Cooperation Programme for 
 2017-2018 under Grant No. 16-1-2016-0027.}

\vspace*{-10mm}

\begin{abstract}
We give sufficient conditions on the offspring, the initial and the immigration distributions under
 which a second-order Galton--Watson process with immigration is regularly varying.
\end{abstract}

\section{Introduction}
\label{section_intro}

Second-order Galton--Watson processes (without immigration) have been introduced by Kashikar and
 Deshmukh \cite{KasDes,KasDes2} and Kashikar \cite{Kas} for modeling the swine flu data for Pune, 
 India and La-Gloria, Mexico. 
They studied their basic probabilistic properties such as a formula for their probability generator
 function, probability of extinction, long run behavior and conditional least squares estimation of
 the offspring means.
 
In a second-order Galton--Watson branching model with immigration, an individual reproduces at age $1$
 and also at age $2$, and then it dies immediately.
For each \ $n \in \NN$, \ the population consists again of the offsprings born at time \ $n$ \ and the
 immigrants arriving at time \ $n$.
\ For each \ $n, i, j \in \NN$, \ the number of offsprings produced at time \ $n$ \ by the
 \ $i^\mathrm{th}$ \ individual of the \ $(n-1)^\mathrm{th}$ \ generation and by the \ $j^\mathrm{th}$
 \ individual of the \ $(n-2)^\mathrm{nd}$ \ generation will be denoted by \ $\xi_{n,i}$ \ and
 \ $\eta_{n,j}$, \ respectively, and \ $\vare_n$ \ denotes the number of immigrants in the
 \ $n^\mathrm{th}$ \ generation.
Then, for the population size \ $X_n$ \ of the \ $n^\mathrm{th}$ \ generation, we have
 \begin{equation}\label{2oGWI}
  X_n = \sum_{i=1}^{X_{n-1}} \xi_{n,i} + \sum_{j=1}^{X_{n-2}} \eta_{n,j} + \vare_n ,
  \qquad n \in \NN ,
 \end{equation}
 where \ $X_{-1}$ \ and \ $X_0$  \ are non-negative integer-valued random variables (the initial
 population sizes).
Here \ $\{X_{-1}, X_0, \, \xi_{n,i}, \, \eta_{n,j}, \, \vare_n : n, i, j \in \NN\}$ \ are supposed to
 be independent non-negative integer-valued random variables, and
 \ $\{\xi_{n,i} : n, i \in \NN\}$, \ $\{\eta_{n,j} : n, j \in \NN\}$ \ and \ $\{\vare_n : n \in \NN\}$
 \ are supposed to consist of identically distributed random variables, respectively.
Note that the number of individuals alive at time \ $n \in \ZZ_+$ \ is \ $X_n + X_{n-1}$, \ which
 can be larger than the population size \ $X_n$ \ of the \ $n^\mathrm{th}$ \ generation, since the
 individuals of the population at time \ $n-1$ \ are still alive at time \ $n$, \ because they can
 reproduce also at age \ $2$.
\ The stochastic process \ $(X_n)_{n\geq -1}$ \ given by \eqref{2oGWI} is called a second-order
 Galton--Watson process with immigration or a Generalized Integer-valued AutoRegressive process
 of order 2 (GINAR(2) process), see, e.g., Latour \cite{Lat}.
Especially, if \ $\xi_{1,1}$ \ and \ $\eta_{1,1}$ \ are Bernoulli distributed random variables, then
 \ $(X_n)_{n\geq -1}$ \ is also called an Integer-valued AutoRegressive process of order 2
 (INAR(2) process), see, e.g., Du and Li \cite{DuLi}.
If \ $\vare_1 = 0$, \ then we say that \ $(X_n)_{n\geq-1}$ \ is a second-order Galton--Watson
 process without immigration.
A process given in \eqref{2oGWI} with the special choice of \ $\eta_{1,1} = 0$ \ gives back a usual
 Galton--Watson process with immigration.

In Section \ref{section_2GWI}, we present conditions for the offspring, the initial and the
 immigration distributions, under which the corresponding second-order Galton--Watson process with
 immigration is regularly varying, extending the results of Barczy et al.\ \cite{BarBosPap} on usual
 Galton--Watson processes with immigration.
Note that these results can be easily generalized with the same techniques to higher order 
 Galton--Watson processes with immigration.

\section{Preliminaries}
\label{section_prel}

Let \ $\ZZ_+$, \ $\NN$, \ $\RR$, \ $\RR_+$ \ and \ $\RR_{++}$ \ denote the set of non-negative
 integers, positive integers, real numbers, non-negative real numbers and positive real numbers,
 respectively.
Equality in distributions of random variables is denoted by \ $\distre$.
\ For functions \ $f : \RR_{++} \to \RR_{++}$ \ and \ $g : \RR_{++} \to \RR_{++}$, \  by the notation
 \ $f(x) \sim g(x)$, \ $f(x) = \oo(g(x))$ \ and \ $f(x) = \OO(g(x))$ \ as \ $x \to \infty$, \ we mean
 that \ $\lim_{x\to\infty} \frac{f(x)}{g(x)} = 1$, \ $\lim_{x\to\infty} \frac{f(x)}{g(x)} = 0$ \ and
 \ $\limsup_{x\to\infty} \frac{f(x)}{g(x)} < \infty$, \ respectively.
Moreover, by the notation \ $f(x) \asymp g(x)$ \ as \ $x \to \infty$, \ we mean that
 \ $f(x) = \OO(g(x))$ \ and \ $g(x) = \OO(f(x))$ \ as \ $x \to \infty$.

First, we recall some results for second-order Galton--Watson processes without or with immigration.
Let \ $(X_n)_{n\geq -1}$ \ be a second-order Galton--Watson process with immigration given in
 \eqref{2oGWI}, and let us introduce the random vectors
 \begin{align}\label{help15}
  \bY_n := \begin{bmatrix}
            Y_{n,1} \\
            Y_{n,2} \\
           \end{bmatrix}
        := \begin{bmatrix}
            X_n \\
            X_{n-1} \\
           \end{bmatrix} ,
  \qquad n \in \ZZ_+ .
 \end{align}
Then we have
 \[
   \bY_n = \sum_{i=1}^{Y_{n-1,1}}
            \begin{bmatrix} \xi_{n,i} \\ 1 \end{bmatrix}
           + \sum_{j=1}^{Y_{n-1,2}}
              \begin{bmatrix} \eta_{n,j} \\ 0 \end{bmatrix}
           + \begin{bmatrix} \vare_n \\ 0 \end{bmatrix} ,
   \qquad n \in \NN ,
 \]
 hence \ $(\bY_n)_{n\in\ZZ_+}$ \ is a (special) 2-type Galton--Watson process with immigration.
In fact, the type  1 and  2 individuals are identified with individuals of age \ $0$ \ and \ $1$,
 \ respectively, and for each \ $n, i, j \in \NN$, \ at time \ $n$, \ the \ $i^\mathrm{th}$
 \ individual of type  1 of the \ $(n-1)^\mathrm{th}$ \ generation produces \ $\xi_{n,i}$
 \ individuals of type  1 and exactly one individual of type 2, and the \ $j^\mathrm{th}$
 \ individual of type 2 of the \ $(n-1)^\mathrm{th}$ \ generation produces \ $\eta_{n,j}$
 \ individuals of type 1 and no individual of type 2.

For notational convenience, let \ $\xi$, \ $\eta$ \ and \ $\vare$ \ be random variables such that
 \ $\xi \distre \xi_{1,1}$, \ $\eta \distre \eta_{1,1}$ \ and \ $\vare \distre \vare_1$.
\ If \ $m_\xi := \EE(\xi) \in \RR_+$, \ $m_\eta := \EE(\eta) \in \RR_+$,
 \ $m_\vare := \EE(\vare) \in \RR_+$, \ $\EE(X_0) \in \RR_+$ \ and \ $\EE(X_{-1}) \in \RR_+$, \ then
 for each \ $n \in \ZZ_+$, \ \eqref{2oGWI} implies
 \[
   \EE(X_n \mid \cF_{n-1}) = X_{n-1} m_\xi + X_{n-2} m_\eta + m_\vare , \qquad  n \in \NN ,
 \]
 where \ $\cF_n := \sigma(X_{-1}, X_0, \dots, X_n)$, \ $n \in \ZZ_+$.
\ Consequently,
 \[
   \EE(X_n) = m_\xi \EE(X_{n-1}) + m_\eta \EE(X_{n-2}) + m_\vare , \qquad n \in \NN ,
 \]
 which can be written in the matrix form
 \begin{equation}\label{recEX_n}
  \begin{bmatrix}
   \EE(X_n) \\
   \EE(X_{n-1})
  \end{bmatrix}
  = \bM_{\xi, \eta}
    \begin{bmatrix}
     \EE(X_{n-1}) \\
     \EE(X_{n-2})
    \end{bmatrix}
    + \begin{bmatrix}
       m_\vare \\
       0
      \end{bmatrix} , \qquad n \in \NN ,
 \end{equation}
 with
 \begin{equation}\label{bA}
  \bM_{\xi, \eta} := \begin{bmatrix}
          m_\xi & m_\eta \\
          1 & 0
         \end{bmatrix} .
 \end{equation}
Note that \ $\bM_{\xi, \eta}^\top$ \ is nothing else but the so-called mean matrix of the 2-type
 Galton--Watson process \ $(\bY_n)_{n\in\ZZ_+}$ \ given in \eqref{help15}.
Thus, we conclude
 \[
   \begin{bmatrix}
    \EE(X_n) \\
    \EE(X_{n-1})
   \end{bmatrix}
   = \bM_{\xi, \eta}^n
     \begin{bmatrix}
      \EE(X_0) \\
       \EE(X_{-1})
     \end{bmatrix}
     + \sum_{k=1}^n
        \bM_{\xi, \eta}^{n-k}
        \begin{bmatrix}
         m_\vare \\
         0
        \end{bmatrix} ,
   \qquad n \in \NN .
 \]

\section{Tail behavior of second-order Galton--Watson processes with immigration}
\label{section_2GWI}

First we consider the case of regularly varying offspring distributions.

\begin{Thm}\label{2GW_xi}
Let \ $(X_n)_{n\geq-1}$ \ be a second-order Galton--Watson process with immigration such \ $\xi$ \ and
 \ $\eta$ \ are regularly varying with index \ $\alpha \in [1, \infty)$,
 \ $\PP(\xi > x) \asymp \PP(\eta > x)$ \ as \ $x\to\infty$, \ $m_\xi, m_\eta \in \RR_{++}$, \ and
 there exists $r \in (\alpha, \infty)$ \ such that \ $\EE(X_0^r)<\infty$, \ $\EE(X_{-1}^r)<\infty$
 \ and \ $\EE(\vare^r)<\infty$. 
\ Suppose that \ $\PP(X_0 = 0) < 1$, \ or \ $\PP(X_{-1} = 0) < 1$ \ or \ $\PP(\vare = 0) < 1$.
\ Then for all \ $n \in \NN$,
 \begin{align*}
  \PP(X_n > x)
  &\sim \begin{bmatrix} \EE(X_0) \\ \EE(X_{-1}) \end{bmatrix}^\top
        \sum_{k=0}^{n-1} m_k^\alpha \bigl(\bM_{\xi, \eta}^{n-k-1}\bigr)^\top
        \begin{bmatrix} \PP(\xi > x) \\ \PP(\eta > x) \end{bmatrix} \\
  &\quad\:
        +
        \begin{bmatrix} m_\vare \\ 0 \end{bmatrix}^\top
        \sum_{i=1}^{n-1} \sum_{j=0} ^{n-i-1} m_j^\alpha \bigl(\bM_{\xi, \eta}^{n-j-1}\bigr)^\top
        \begin{bmatrix} \PP(\xi > x) \\ \PP(\eta > x) \end{bmatrix}
   \qquad \text{as $x\to\infty$,}
 \end{align*}
 where \ $m_0 := 1$ \ and
 \begin{align}\label{m_n}
  m_k := \frac{\lambda_+^{k+1}-\lambda_-^{k+1}}{\lambda_+-\lambda_-} , \qquad
  \lambda_+ := \frac{m_\xi+\sqrt{m_\xi^2+4m_\eta}}{2} , \qquad
  \lambda_- := \frac{m_\xi-\sqrt{m_\xi^2+4m_\eta}}{2}
 \end{align}
 for \ $k \in \NN$.
\ As a consequence, \ $X_n$ \ is regularly varying with index \ $\alpha$ \ for all \ $n \in \NN$.
\end{Thm}

Note that \ $\lambda_+$ \ and \ $\lambda_-$ \ are the eigenvalues of the matrix \ $\bM_{\xi, \eta}$
 \ given in \eqref{bA} related to the recursive formula \eqref{recEX_n} for the expectations
 \ $\EE(X_n)$, \ $n \in \NN$.
\ One can use Lemma B1 in Barczy et al.\ \cite{BarBosPap2} for evaluating \ $\bM_{\xi, \eta}^\ell$,
 \ $\ell \in \ZZ_+$.
\ Moreover, for all \ $k \in \ZZ_+$, \ the assumptions \ $m_\xi \in \RR_{++}$ \ and
 \ $m_\eta \in \RR_{++}$ \ imply \ $m_k \in \RR_{++}$.
\ Further, by Lemma B1 in Barczy et al.\ \cite{BarBosPap2}, for all \ $k \in \ZZ_+$, \ we have
 \ $m_k = \EE(V_{k,0})$, \ where \ $(V_{n,0})_{n\geq-1}$ \ is a second-order Galton--Watson process
 without immigration with initial values \ $V_{0,0} = 1$, \ $V_{-1,0} = 0$ \ and with the same
 offspring distributions as \ $(X_n)_{n\geq-1}$.

\noindent{\bf Proof of Theorem \ref{2GW_xi}.}
First, observe that for each \ $n \in \NN$, \ by the additive property \eqref{2GWI_additive}, we have
 \[
   X_n = V_0^{(n)}(X_0, X_{-1}) + \sum_{i=1}^n V_i^{(n-i)}(\vare_i, 0),
 \]
 where we recall that \ $V_0^{(n)}(X_0, X_{-1})$ \ represents the number of newborns at time \ $n$,
 \ resulting from the initial individuals \ $X_0$ \ at time \ $0$, \ and for each
 \ $i \in \{1, \ldots, n\}$, \ $V_i^{(n-i)}(\vare_i, 0)$ \ represents the number of newborns at time
 \ $n$, \ resulting from the immigration \ $\vare_i$ \ at time \ $i$.
\ Moreover, here, by \eqref{2GW_additive}, we have that
 \begin{equation}\label{V0GWI}
  V_0^{(n)}(X_0, X_{-1})
  \distre \sum_{i=1}^{X_0} \zeta_{i,0}^{(n)} + \sum_{j=1}^{X_{-1}} \zeta_{j,-1}^{(n)} ,
 \end{equation}
 where \ $\{\zeta_{i,0}^{(n)} : i \in \NN\}$ \ are independent copies of \ $V_{n,0}$ \ and
 \ $\{\zeta_{j,-1}^{(n)} : j \in \NN\}$ \ are independent copies of \ $V_{n,-1}$, \ where
 \ $(V_{k,0})_{k\geq-1}$ \ and \ $(V_{k,-1})_{k\geq-1}$ \ are second-order Galton--Watson processes
 (without immigration) with initial values \ $V_{0,0} = 1$, \ $V_{-1,0} = 0$, \ $V_{0,-1} = 0$ \ and
 \ $V_{-1,-1} = 1$, \ and with the same offspring distributions as \ $(X_k)_{k\geq-1}$.

First we prove that for all $n \in \NN$
 \begin{equation}\label{VV}
  \begin{bmatrix} \PP(V_{n,0} > x) \\ \PP(V_{n,-1} > x) \end{bmatrix}
  \sim \sum_{k=0}^{n-1} m_k^\alpha \bigl(\bM_{\xi, \eta}^{n-k-1}\bigr)^\top
       \begin{bmatrix} \PP(\xi > x) \\ \PP(\eta > x) \end{bmatrix}
  \qquad \text{as \ $x \to \infty$,}
 \end{equation}
 understanding the relation \ $\sim$ \ as \ $x \to \infty$ \ for both coordinates.
We proceed by induction on \ $n$.
\ For \ $n = 1$, \ the statement readily follows, since \ $V_{1,0} = \xi_{1,1} \distre \xi$ \ and
 \ $V_{1,-1} = \eta_{1,1} \distre \eta$.
\ Now let us assume that \eqref{VV} hold for \ $1, \ldots, n - 1$, \ where \ $n \geq 2$.
\ Since
 \[
   \Biggl(\begin{bmatrix} V_{k,0} \\ V_{k-1,0} \end{bmatrix}\Biggr)_{k\in\NN} , \qquad
   \Biggl(\begin{bmatrix} V_{k,-1} \\ V_{k-1,-1} \end{bmatrix}\Biggr)_{k\in\NN}
 \]
 are time homogeneous Markov processes with initial vectors
 \[
   \begin{bmatrix} V_{1,0} \\ V_{0,0} \end{bmatrix} = \begin{bmatrix} \xi_{1,1} \\ 1 \end{bmatrix} ,
   \qquad
   \begin{bmatrix} V_{1,-1} \\ V_{0,-1} \end{bmatrix}
   = \begin{bmatrix} \eta_{1,1} \\ 0 \end{bmatrix} ,
 \]
 we have
 \[
   \begin{bmatrix} V_{k,0} \\ V_{k-1,0} \end{bmatrix}
   \distre \begin{bmatrix} V^{(k-1)}(\xi_{1,1}, 1) \\ V^{(k-2)}(\xi_{1,1}, 1) \end{bmatrix} , \qquad
   \begin{bmatrix} V_{k,-1} \\ V_{k-1,-1} \end{bmatrix}
   \distre \begin{bmatrix} V^{(k-1)}(\eta_{1,1}, 0) \\ V^{(k-2)}(\eta_{1,1}, 0) \end{bmatrix} ,
 \]
 where \ $(V^{(k)}(\xi_{1,1}, 1))_{k\in\ZZ_+}$ \ and \ $(V^{(k)}(\eta_{1,1}, 0))_{k\in\ZZ_+}$ \ are
 (second-order) Galton--Watson processes (without immigration) with initial values
 \ $V^{(0)}(\xi_{1,1}, 1) = \xi_{1,1}$, \ $V^{(-1)}(\xi_{1,1}, 1) = 1$,
 \ $V^{(0)}(\eta_{1,1}, 0) = \eta_{1,1}$ \ and \ $V^{(-1)}(\eta_{1,1}, 0) = 0$.
\ Applying again the additive property \eqref{2GW_additive}, we obtain
 \[
   V_{n,0} \distre V^{(n-1)}(\xi_{1,1}, 1)
   \distre \sum_{i=1}^{\xi_{1,1}} \zeta_{i,0}^{(n-1)} + \zeta_{1,-1}^{(n-1)} , \qquad
   V_{n,-1} \distre V^{(n-1)}(\eta_{1,1}, 0)
   \distre \sum_{i=1}^{\eta_{1,1}} \zeta_{i,0}^{(n-1)} ,
 \]
 where \ $\{\zeta_{i,0}^{(n-1)} : i \in \NN\}$ \ and \ $\zeta_{1,-1}^{(n-1)}$ \ are independent
 copies of \ $V_{n-1,0}$ \ and \ $V_{n-1,-1}$, \ respectively, such that
 \ $\{\xi_{1,1}, \eta_{1,1}, \zeta_{i,0}^{(n-1)}, \zeta_{1,-1}^{(n-1)} : i \in \NN\}$ \ are
 independent.
First note that \ $\PP(\zeta_{1,0}^{(n-1)} > x) = \OO(\PP(\xi > x))$ \ and
 \ $\PP(\zeta_{1,0}^{(n-1)} > x) = \OO(\PP(\eta > x))$ \ as \ $x \to \infty$.
\ Indeed, using the induction hypothesis and the condition \ $\PP(\xi > x) \asymp \PP(\eta > x)$ \ as
 \ $x\to\infty$, \ we obtain
 \begin{align*}
  \limsup_{x\to\infty} \frac{\PP(\zeta_{1,0}^{(n-1)} > x)}{\PP(\xi > x)}
  &= \limsup_{x\to\infty} \frac{\PP(V_{n-1,0} > x)}{\PP(\xi > x)} \\
  &= \limsup_{x\to\infty}
      \begin{bmatrix} 1 \\ 0 \end{bmatrix}^\top
      \sum_{k=0}^{n-2} m_k^\alpha \bigl(\bM_{\xi, \eta}^{n-k-2}\bigr)^\top
      \begin{bmatrix} 1 \\ \frac{\PP(\eta > x)}{\PP(\xi > x)} \end{bmatrix}
   < \infty .
 \end{align*}
The statement \ $\PP(\zeta_{1,0}^{(n-1)} > x) = \OO(\PP(\eta > x))$ \ as \ $x \to \infty$ \ can be
 proven similarly.
Now we can apply Proposition \ref{DFK_corr}, and we obtain
 \begin{align*}
  \PP\left(\sum_{i=1}^{\xi_{1,1}} \zeta_{i,0}^{(n-1)} > x\right)
  &\sim \EE(\xi_{1,1}) \PP(\zeta_{1,0}^{(n-1)} > x)
        + \PP\left(\xi_{1,1} > \frac{x}{\EE(\zeta_{1,0}^{(n-1)})}\right) \\
  &\sim m_\xi \PP(V_{n-1,0} > x) + m_{n-1}^\alpha \PP(\xi > x)
 \end{align*}
 as \ $x \to \infty$, \ yielding that \ $\sum_{i=1}^{\xi_{1,1}} \zeta_{i,0}^{(n-1)}$ \ is regularly
 varying with index \ $\alpha$, \ and
 \begin{align*}
  \PP\left(\sum_{i=1}^{\eta_{1,1}} \zeta_{i,0}^{(n-1)} > x\right)
  &\sim \EE(\eta_{1,1}) \PP(\zeta_{1,0}^{(n-1)} > x)
        + \PP\left(\eta_{1,1} > \frac{x}{\EE(\zeta_{1,0}^{(n-1)})}\right) \\
  &\sim m_\eta \PP(V_{n-1,0} > x) + m_{n-1}^\alpha \PP(\eta > x)
 \end{align*}
 as \ $x \to \infty$, \ yielding that \ $\sum_{i=1}^{\eta_{1,1}} \zeta_{i,0}^{(n-1)}$ \ is regularly
 varying with index \ $\alpha$.
\ Consequently, by the convolution property, we get
 \[
   \sum_{i=1}^{\xi_{1,1}} \zeta_{i,0}^{(n-1)} + \zeta_{1,-1}^{(n-1)}
   \sim m_\xi \PP(V_{n-1,0} > x) + m_{n-1}^\alpha \PP(\xi > x) + \PP(V_{n-1,-1} > x)
 \]
 as \ $x \to \infty$.
\ Writing these relations in a matrix form, we conclude
 \[
   \begin{bmatrix} \PP(V_{n,0} > x) \\ \PP(V_{n,-1} > x) \end{bmatrix}
   \sim \begin{bmatrix} m_\xi & 1 \\ m_\eta & 0 \end{bmatrix}
        \begin{bmatrix} \PP(V_{n-1,0} > x) \\ \PP(V_{n-1,-1} > x) \end{bmatrix}
        + m_{n-1}^\alpha \begin{bmatrix} \PP(\xi > x) \\ \PP(\eta > x) \end{bmatrix}
  \qquad \text{as \ $x \to \infty$.}
 \]
Using the induction hypothesis and that \ $V_{1,0} \distre \xi$ \ and
 \ $V_{1,-1} \distre \eta$, \ we obtain \eqref{VV}. 
From this we can conclude that \ $V_{n,0}$ \ and \ $V_{n,-1}$ \ are regularly varying with index
 \ $\alpha$, \ and now we can apply Proposition \ref{FGAMSRS}, and we have that 
 \begin{equation}\label{RVV0N}
  \PP(V_0^{(n)}(X_0, X_{-1}) > x)
  \sim \begin{bmatrix} \EE(X_0) \\ \EE(X_{-1}) \end{bmatrix}^\top
       \sum_{k=0}^{n-1} m_k^\alpha \bigl(\bM_{\xi, \eta}^{n-k-1}\bigr)^\top
       \begin{bmatrix} \PP(\xi > x) \\ \PP(\eta > x) \end{bmatrix} 
 \end{equation}
 as \ $x \to \infty$.
\ Now note that \ $V_i^{(n-i)}(\vare_i, 0) \distre \sum_{j=1}^{\vare_i} \zeta_{j,0}^{(n-i)}$ \ for all
 \ $i \in \{1, \ldots, n\}$. 
\ For \ $i = n$, \ we have \ $V_n^{(0)}(\vare_n, 0) = \vare_n$, \ and with the moment assumption on
 \ $\vare$, \ by Lemma \ref{exposv}, we get \ $\PP(\vare > x) = \oo(\PP(\xi>x))$ \ and
 \ $\PP(\vare > x) = \oo(\PP(\eta > x))$ \ as \ $x \to \infty$. 
\ Otherwise, we can apply \eqref{RVV0N}, and we have that for all \ $i \in \{1, \ldots, n - 1\}$, 
 \begin{equation}\label{RVEPS}
  \PP(V_i^{(n-i)}(\vare_i, 0) > x)
  \sim \begin{bmatrix} m_\vare \\ 0 \end{bmatrix}^\top
       \sum_{j=0} ^{n-i-1} m_j^\alpha \bigl(\bM_{\xi, \eta}^{n-j-1}\bigr)^\top
       \begin{bmatrix} \PP(\xi > x) \\ \PP(\eta > x) \end{bmatrix} ,
 \end{equation}
 hence \ $V_i^{(n-i)}(\vare_i, 0)$ \ is regularly varying with index \ $\alpha$. 
\ Now, using the convolution property in Lemma \ref{Lem_conv}, we obtain the desired tail behavior for
 \ $X_n$.
\proofend
		
\begin{Rem} \label{Rem1}
Similar statements can be derived with a similar proof, when we assume that only one offspring
 distribution is regularly varying with index \ $\alpha$, \ and for the other we suppose a finite
 $r$-th moment with \ $r \in (\alpha, \infty)$. 
\ In this case, in the tail behavior the vector \ $(\PP(\xi > x), 0)^\top$ \ will appear instead of
 \ $(\PP(\xi > x), \PP(\eta > x))^\top$, \ everything else will remain the same.  
To prove this proposition, we can follow the lines of the proof of Theorem \ref{2GW_xi}, there is only
 one modification that we have to make: when determining the asymptotics of the probability
 \ $\PP\bigl(\sum_{i=1}^{\eta_{1,1}} \zeta_{i,0}^{(n-1)} > x\bigr)$, \ we have to apply Lemma
 \ref{DFK_corr0} instead of Lemma \ref{DFK_corr}.

Similarly, if we suppose that \ $\eta$ \ is regularly varying with index \ $\alpha \in[1, \infty)$,
 \ $m_\eta \in \RR_{++}$ \ and there exits \ $r \in (\alpha, \infty)$ \ such that
 \ $\EE(X_0^r) < \infty$, \ $\EE(X_{-1}^r) < \infty$, \ $\EE(\xi^r) < \infty$ \ and
 \ $\EE(\vare^r) < \infty$, \ then in the tail behavior the vector \ $(0, \PP(\eta > x))^\top$ \ will
 appear instead of \ $(\PP(\xi > x), 0)^\top$, \ everything else will remain the same.
\proofend
\end{Rem}

Next, we consider the case of regularly varying initial distributions.

\begin{Thm}\label{2GW_X_0_X_-1}
Let \ $(X_n)_{n\geq-1}$ \ be a second-order Galton--Watson process with immigration such that \ $X_0$
 \ and \ $X_{-1}$ \ are regularly varying with index \ $\beta \in \RR_+$, \ and there exists
 \ $r \in (1 \lor \beta, \infty)$ \ with \ $\EE(\xi^r) < \infty$, \ $\EE(\eta^r) < \infty$ \ and
 \ $\EE(\vare^r) < \infty$.
\ Suppose that \ $\PP(\xi = 0) < 1$ \ or \ $\PP(\eta = 0) < 1$.
\ Then for all \ $n \in \NN$,
 \[
   \PP(X_n > x) \sim m_n^\beta \PP(X_0 > x) + m_{n-1}^\beta m_\eta^\beta \PP(X_{-1} > x)
   \qquad \text{as \ $x \to \infty$,}
 \]
 where \ $m_n$, \ $n \in \NN$, \ are given in \eqref{m_n}, yielding that \ $X_n$ \ is regularly
 varying with index \ $ \beta$.
\end{Thm}

\noindent{\bf Proof.}
Let us fix \ $n \in \NN$.
\ In view of the additive properties \eqref{2GWI_additive} and \eqref{V0GWI}, and the convolution
 property of regularly varying functions described in Lemma \ref{Lem_conv}, it is sufficient to prove
 \begin{equation}\label{ac}
  \PP\Biggl(\sum_{i=1}^{X_0} \zeta_{i,0}^{(n)} > x\Biggr) \sim m_n^{\beta} \PP(X_0 > x) , \qquad
  \PP\Biggl(\sum_{j=1}^{X_{-1}} \zeta_{j,-1}^{(n)} > x\Biggr)
  \sim m_{n-1}^{\beta} m_\eta^{\beta} \PP(X_{-1} > x)
 \end{equation}
 as \ $x \to \infty$, \ because by Lemma \ref{Lem_seged_momentr},
 \ $\EE\bigl(V_i^{(n-i)}(\vare_i, 0)^r\bigr)<\infty$ \ for all \ $i \in \{1, \ldots, n\}$. 
\ Also from Lemma \ref{Lem_seged_momentr}, we can conclude that
 \ $\EE((\zeta_{1,0}^{(n)})^r) < \infty$ \ and \ $\EE((\zeta_{1,-1}^{(n)})^r)<\infty$.
\ Then the desired relations in \eqref{ac} follow from Proposition \ref{FGAMSRS}, since
 \ $\EE(\zeta_{1,0}^{(n)}) = m_n$, \ $\EE(\zeta_{1,-1}^{(n)}) = m_{n-1} m_\eta \in\RR_{++}$,
 \ $n \in \NN$, \ by Lemma B1 in Barczy et al.\ \cite{BarBosPap2}.
\proofend
	
\begin{Rem} \label{Rem2}
We can prove similar statement if we only assume that \ $X_0$ \ is regularly varying with index
 \ $\beta \in \RR_+$, \ and we make some moment assumption on \ $X_{-1}$, \ then only the term
 \ $m_n^\beta \PP(X_0 > x)$ \ will appear in the tail behavior.
To prove this, we can follow the lines of the proof of Theorem \ref{2GW_X_0_X_-1}, we just need to
 observe that now, by Lemma \ref{Lem_seged_momentr}, we have that
 \ $\EE\bigl(\bigl(\sum_{j=1}^{X_{-1}} \zeta_{j,-1}^{(n)}\bigr)^r\bigr) < \infty$, \ and consequently,
 by Lemma \ref{Lem_conv}, this term will not affect the tail behavior of the process.

Also, if we change here the assumptions on \ $X_0$ \ and \ $X_{-1}$ \ with each other (with leaving
 the other assumptions the same), it can be proven similarly that in this case,
 \ $\PP(X_n > x) \sim m_{n-1}^\beta m_\eta^\beta \PP(X_{-1} > x)$ \ as \ $x \to \infty$.
\proofend
\end{Rem}

With the combination of the previous theorems, we obtain the following result.

\begin{Thm}\label{2GW_X_0_X_-1_xi_3}
Let \ $(X_n)_{n\geq-1}$ \ be a second-order Galton--Watson process with immigration such that \ $X_0$,
 \ $X_{-1}$, \ $\xi$ \ and \ $\eta$ \ are regularly varying with index \ $\alpha \in [1, \infty)$,
 \ $\EE(X_0), \EE(X_{-1}), m_\xi, m_\eta \in \RR_{++}$, \ $\PP(\xi > x)  = \OO(\PP(X_0 > x))$,
 \ $\PP(\xi > x) = \OO(\PP(X_{-1} > x))$, \ $\PP(\eta > x)  = \OO(\PP(X_0 > x))$,
 \ $\PP(\eta > x) = \OO(\PP(X_{-1} > x))$ \ and \ $\PP(\xi > x) \asymp \PP(\eta > x)$ \ as
 \ $\to \infty$, \ and there exists \ $r \in (\alpha, \infty)$ \ such that \ $\EE(\vare^r) < \infty$. 
\ Then for all \ $n \in \NN$,
 \begin{align*}
  \PP(X_n > x)
  &\sim \begin{bmatrix} \EE(X_0) \\ \EE(X_{-1}) \end{bmatrix}^\top
        \sum_{k=0}^{n-1} m_k^\alpha \bigl(\bM_{\xi, \eta}^{n-k-1}\bigr)^\top
        \begin{bmatrix} \PP(\xi > x) \\ \PP(\eta > x) \end{bmatrix} \\  
  &\quad  +
  \begin{bmatrix} m_\vare \\ 0 \end{bmatrix}^\top
  \sum_{i=1}^{n-1} \sum_{j=0} ^{n-i-1} m_j^\alpha \bigl(\bM_{\xi, \eta}^{n-j-1}\bigr)^\top
  \begin{bmatrix} \PP(\xi > x) \\ \PP(\eta > x) \end{bmatrix} \\  
  &\quad
        + m_n^\alpha \PP(X_0 > x) + m_{n-1}^\alpha m_\eta^\alpha \PP(X_{-1} > x)
   \qquad \text{as \ $x \to \infty$,}
 \end{align*}
 where \ $m_n$, \ $n \in \NN$, \ are given in \eqref{m_n}, yielding that \ $X_n$ \ is regularly
 varying with index \ $\alpha$.
\end{Thm}

\noindent{\bf Proof.}
By the additive properties \eqref{2GW_additive} and \eqref{2GWI_additive}, we have
 \[
   X_n \distre \sum_{i=1}^{X_0} \zeta_{i,0}^{(n)} + \sum_{j=1}^{X_{-1}} \zeta_{j,-1}^{(n)}
               + \sum_{i=1}^n V_i^{(n-i)} (\vare_i, 0) .
 \]
From the assumptions it follows that \ $\PP(\zeta_{1,0}^{(n)} > x) = \OO(\PP(X_0 >x))$ \ and
 \ $\PP(\zeta_{1,-1}^{(n)} > x) = \OO(\PP(X_{-1} > x))$ \ as \ $x \to \infty$, \ so by Lemma
 \ref{DFK_corr}, we get 
 \[
   \PP\Biggl(\sum_{i=1}^{X_0} \zeta_{i,0}^{(n)} > x\Biggr)
   \sim \EE(X_0)\PP(V_{n,0} > x) + m_n^\alpha \PP(X_0 > x) \qquad \text{as \ $x \to \infty$}
 \]
 and
 \[
   \PP\Biggl(\sum_{j=1}^{X_{-1}} \zeta_{j, -1}^{(n)} > x\Biggr)
   \sim \EE(X_{-1}) \PP(V_{n,-1} > x) + m_{n-1}^\alpha m_\eta^\alpha \PP(X_{-1} > x)
   \qquad \text{as \ $x \to \infty$.}
 \]
From \eqref{VV} we know the asymptotics of the probabilities \ $\PP(V_{n,0} > x)$ \ and
 \ $\PP(V_{n,-1} > x)$ \ as \ $x \to \infty$, \ and from \eqref{RVEPS}, we know the asymptotics of the
 probabilities \ $\PP(V_i^{(n-i)}(\vare_i, 0) > x)$ \ as \ $x \to \infty$ \ for each
 \ $i \in \{1, \ldots, n - 1\}$, \ and for \ $i = n$, \ we have \ $V_n^{(0)}(\vare_n,0) = \vare_n$. 
\ Applying the convolution property in Lemma \ref{Lem_conv}, we conclude the statement.
\proofend

\begin{Rem} \label{Rem3}
We can relax some conditions on the regular variation of some of the random variables appearing in 
 Theorem  \ref{2GW_X_0_X_-1_xi_3} and replace them with moment assumptions, as described in Remark
 \ref{Rem1} and Remark \ref{Rem2}. 
The tail behavior changes accordingly to these remarks.
\proofend
\end{Rem}

Now we consider the case of regularly varying immigration.

\begin{Thm} \label{2GWI_vare}
Let \ $(X_n)_{n\in\ZZ_+}$ \ be a second-order Galton--Watson process with immigration such that
 \ $\vare$ \ is regularly varying with index \ $\gamma \in \RR_+$, \ and there exists
 \ $r \in (1 \lor \gamma, \infty)$ \ with \ $\EE(\xi^r) < \infty$, \ $\EE(\eta^r) < \infty$,
 \ $\EE(X_0^r)<\infty$ \ and \ $\EE(X_0^r) < \infty$.
\ Suppose that \ $\PP(\xi = 0) < 1$ \ or \ $\PP(\eta = 0) < 1$.
\ Then for each \ $n \in \NN$, \ we have
 \[
   \PP(X_n > x) \sim \sum_{i=1}^n m_{n-i}^\gamma \PP(\vare > x) \qquad \text{as \ $x \to \infty$,}
 \]
 and hence \ $X_n$ \ is also regularly varying with index \ $\gamma$.
\end{Thm}

\noindent{\bf Proof.}
We use the representation \eqref{2GWI_additive}. 
By Lemma \ref{Lem_seged_momentr}, we have that \ $\EE(V_0^{(n)}(X_0, X_{-1})^r) < \infty$. 
\ For each \ $i \in \{1, \ldots, n\}$, \ by Lemma \ref{DFK_corr0}, we have
 \ $\PP(V_i^{(n-i)}(\vare_i, 0) > x) \sim m_{n-i}^\gamma \PP(\vare > x)$ \ as \ $x \to \infty$.
\ Applying the convolution property in Lemma \ref{Lem_conv}, we conclude the statement.
\proofend

\begin{Thm}\label{2GWIXIEPS}
Let \ $(X_n)_{n\geq-1}$ \ be a second-order Galton--Watson process with immigration such \ $\xi$,
 \ $\eta$ \ and \ $\vare$ \ are regularly varying with index \ $\alpha \in [1, \infty)$,
 \ $\PP(\xi > x) \asymp \PP(\eta > x)$, \ $\PP(\xi > x) = \OO(\PP(\vare > x))$ \ and
 \ $\PP(\eta > x) = \OO(\PP(\vare>x))$ \ as \ $x \to \infty$, \ $m_\xi, m_\eta, m_\vare \in \RR_{++}$,
 \ and there exists \ $r \in (\alpha, \infty)$ \ such that \ $\EE(X_0^r) < \infty$ \ and
 \ $\EE(X_{-1}^r) < \infty$. 
\ Then for all \ $n \in \NN$,
 \begin{align*}
  \PP(X_n > x)
  &\sim \begin{bmatrix} \EE(X_0) \\ \EE(X_{-1}) \end{bmatrix}^\top
        \sum_{k=0}^{n-1} m_k^\alpha \bigl(\bM_{\xi, \eta}^{n-k-1}\bigr)^\top
        \begin{bmatrix} \PP(\xi > x) \\ \PP(\eta > x) \end{bmatrix} \\
  &\quad\:
        + \begin{bmatrix} m_\vare \\ 0 \end{bmatrix}^\top
          \sum_{i=1}^{n-1} \sum_{j=0} ^{n-i-1} m_j^\alpha \bigl(\bM_{\xi, \eta}^{n-j-1}\bigr)^\top
           \begin{bmatrix} \PP(\xi > x) \\ \PP(\eta > x) \end{bmatrix}
        + \sum_{i=1}^n m_{n-i}^\alpha \PP(\vare>x)
 \end{align*} 
 as \ $x \to \infty$, \ yielding that \ $X_n$ \ is also regularly varying with index \ $\alpha$.
\end{Thm}

\noindent{\bf Proof.}
We use the representation \eqref{2GWI_additive}. 
By the proof of Theorem \ref{2GW_xi} we have that 
 \[
   \PP(V_0^{(n)}(X_0, X_{-1}) > x) \sim \begin{bmatrix} \EE(X_0) \\ \EE(X_{-1}) \end{bmatrix}^\top
   \sum_{k=0}^{n-1} m_k^\alpha \bigl(\bM_{\xi, \eta}^{n-k-1}\bigr)^\top
   \begin{bmatrix} \PP(\xi > x) \\ \PP(\eta > x) \end{bmatrix}
 \]
 as \ $x \to \infty$.
\ By Lemma \ref{DFK_corr}, we can conclude that for all \ $i \in \{1, \ldots, n - 1\}$,
 \[
   \PP(V_i^{(n-i)}(\vare_i, 0) > x) \sim m_\vare \PP(V_{n-i,0} > x) + m_{n-i}^\alpha \PP(\vare > x)
 \]
 as \ $x \to \infty$, \ and for \ $i = n$, \ we have \ $V_n^{(0)}(\vare_n,0) = \vare_n$.
\ Again by the proof of Theorem \ref{2GW_xi}, we have that for all \ $i \{1, \ldots, n - 1\}$,
 \[
   m_\vare \PP(V_{n-i,0} >x)
   \sim \sum_{j=0}^{n-i-1}
         m_j^\alpha \bigl(\bM_{\xi, \eta}^{n-j-1}\bigr)^\top
         \begin{bmatrix} \PP(\xi > x) \\ \PP(\eta > x) \end{bmatrix}
 \]
 as \ $x \to \infty$.
\ Applying the convolution property in Lemma \ref{Lem_conv}, we conclude the statement.
\proofend

\begin{Rem}
Similarly to Remark \ref{Rem1}, we can assume that only one of the offspring distributions is 
 regularly varying, and we suppose that the other satisfies a moment condition. 
The tail behavior changes accordingly to Remark \ref{Rem1}.
\proofend
\end{Rem}

\begin{Thm} \label{2GWIX0VARE}
Let \ $(X_n)_{n\geq-1}$ \ be a second-order Galton--Watson process with immigration such that \ $X_0$,
 \ $X_{-1}$ \ and \ $\vare$ \ are regularly varying with index \ $\beta \in \RR_+$, \ and there
 exists \ $r \in (1 \lor \beta, \infty)$ \ with \ $\EE(\xi^r) < \infty$ \ and
 \ $\EE(\eta^r) < \infty$.
\ Then for all \ $n \in \NN$,
 \[
   \PP(X_n > x)
   \sim m_n^\beta \PP(X_0 > x) + m_{n-1}^\beta m_\eta^\beta \PP(X_{-1} > x)
        + \sum_{i=1}^n m_{n-i}^\beta \PP(\vare > x)
 \]
 as \ $x \to \infty$, \ yielding that \ $X_n$ \ is also regularly varying with index \ $\beta$.
\end{Thm}

\noindent{\bf Proof.}
We use the representation \eqref{2GWI_additive}. By Theorem \ref{2GW_X_0_X_-1} we have that 
 \[
   \PP(V_0^{(n)}(X_0, X_{-1}) > x)
   \sim m_n^\beta \PP(X_0 > x) + m_{n-1}^\beta m_\eta^\beta \PP(X_{-1} > x)
 \]
 as \ $x \to \infty$. 
\ By  the proof of Theorem \ref{2GWI_vare} we know that for each \ $i \in \{1, \ldots, n\}$,  
 \[
   \PP(V_i^{(n-i)}(\vare_i, 0) >x) \sim m_{n-i}^\beta \PP(\vare > x), \qquad \text{as $x\to\infty$.}
 \]
Applying the convolution property in Lemma \ref{Lem_conv}, we conclude the statement.
\proofend

\begin{Rem}
Similarly to Remark \ref{Rem2}, we can assume that only one of the initial distributions are regularly
 varying, and for the other a moment assumption holds. 
The tail behavior changes accordingly to Remark \ref{Rem2}. 
\end{Rem}

\begin{Thm}\label{2GWIX0XIVARE}
If \ $(X_n)_{n\geq-1}$ \ is a second-order Galton--Watson process with immigration such that
 \ $X_0$, \ $X_{-1}$, \ $\xi$, \ $\eta$ \ and \ $\vare$ \ are regularly varying with index
 \ $\alpha \in [1, \infty)$, \ $\EE(X_0), \EE(X_{-1}), m_\xi, m_\vare\in \RR_{++}$,
 \ $m_\eta \in \RR_+$, \ $\PP(\xi > x)  = \OO(\PP(X_0 > x))$, \ $\PP(\xi > x) = \OO(\PP(X_{-1} > x))$,
 \ $\PP(\eta > x) = \OO(\PP(X_0 > x))$, \ $\PP(\eta > x) = \OO(\PP(X_{-1} > x))$,
 \ $\PP(\xi > x) = \OO(\PP(\vare > x))$, \ $\PP(\eta > x) = \OO(\PP(\vare > x))$ \ and
 \ $\PP(\xi > x) \asymp \PP(\eta > x)$ \ as \ $x \to \infty$. 
\ Then for all \ $n \in \NN$,
 \begin{align*}
  \PP(X_n > x)
  &\sim \begin{bmatrix} \EE(X_0) \\ \EE(X_{-1}) \end{bmatrix}^\top
        \sum_{k=0}^{n-1} m_k^\alpha \bigl(\bM_{\xi, \eta}^{n-k-1}\bigr)^\top
        \begin{bmatrix} \PP(\xi > x) \\ \PP(\eta > x) \end{bmatrix} \\  
  &\quad
        + \begin{bmatrix} m_\vare \\ 0 \end{bmatrix}^\top
          \sum_{i=1}^{n-1} \sum_{j=0} ^{n-i-1} m_j^\alpha \bigl(\bM_{\xi, \eta}^{n-j-1}\bigr)^\top
          \begin{bmatrix} \PP(\xi > x) \\ \PP(\eta > x) \end{bmatrix} \\  
  &\quad
        + m_n^\beta \PP(X_0 > x) + m_{n-1}^\alpha m_\eta^\alpha \PP(X_{-1} > x) 
        + \sum_{i=1}^n m_{n-i}^\alpha \PP(\vare>x)
 \end{align*}
 as \ $x \to \infty$, \ yielding that \ $X_n$ \ is also regularly varying with index \ $\alpha$.
\end{Thm}

\noindent{\bf Proof.}
We use the additive property \eqref{2GWI_additive}. 
By applying Theorem \ref{2GW_X_0_X_-1_xi_3} with \ $\vare = 0$, \ we get 
 \begin{align*}
  \PP(V_0^{(n)}(X_0, X_{-1}) > x)
  &\sim \begin{bmatrix} \EE(X_0) \\ \EE(X_{-1}) \end{bmatrix}^\top
        \sum_{k=0}^{n-1} m_k^\alpha \bigl(\bM_{\xi, \eta}^{n-k-1}\bigr)^\top
        \begin{bmatrix} \PP(\xi > x) \\ \PP(\eta > x) \end{bmatrix} \\  
  &\quad
        + m_n^\alpha \PP(X_0 > x) + m_{n-1}^\beta m_\eta^\alpha \PP(X_{-1} > x)
 \end{align*}
 as \ $x\to\infty$. 
\ By Lemma \ref{DFK_corr}, we have that for all \ $i \in \{1, \ldots, n - 1\}$ \ (as it can be seen in
 the proof of Theorem \ref{2GWIXIEPS}),
 \[
   \PP(V_i^{(n-i)(\vare_i, 0)} > x)
   \sim \begin{bmatrix} m_\vare \\ 0 \end{bmatrix}^\top
        \sum_{j=0} ^{n-i-1}
         m_j^\alpha \bigl(\bM_{\xi, \eta}^{n-j-1}\bigr)^\top
         \begin{bmatrix} \PP(\xi > x) \\ \PP(\eta > x) \end{bmatrix} \\  
        + m_{n-i}^\alpha \PP(\vare > x)
 \]
 as \ $x \to \infty$, \ and for \ $i = n$, \ we have that \ $V_n^{(0)} (\vare_n, 0) = \vare_n$.
\ Applying the convolution property in Lemma \ref{Lem_conv}, we conclude the statement.
\proofend

\begin{Rem}
We can also relax the condition of regular variation of one of the offspring distribution and one of
 the initial distribution. 
The tail behavior will change according to Remark \ref{Rem1} and \ref{Rem2}.
\proofend
\end{Rem}

\vspace*{2mm}

\appendix

\noindent{\bf\Large Appendices}

\vspace*{-5mm}

\section{Representations of second-order Galton--Watson processes}
\label{App0}

Note that, for a second-order Galton--Watson process \ $(X_n)_{n\geq -1}$ \ (without immigration), the
 additive (or branching) property of a 2-type Galton--Watson process (without immigration), see, e.g.\
 in Athreya and Ney \cite[Chapter V, Section 1]{AthNey}, together with the law of total probability,
 imply
 \begin{equation}\label{2GW_additive}
  X_n \distre \sum_{i=1}^{X_0} \zeta_{i,0}^{(n)} + \sum_{j=1}^{X_{-1}} \zeta_{j,-1}^{(n)} ,
 \end{equation}
 where \ $\bigl\{X_0, X_{-1}, \zeta_{i,0}^{(n)}, \zeta_{j,-1}^{(n)} : i, j \in \NN\bigr\}$ \ are
 independent random variables such that \ $\{\zeta_{i,0}^{(n)} : i \in \NN\}$ \ are independent copies
 of \ $V_{n,0}$ \ and \ $\{\zeta_{j,-1}^{(n)} : j \in \NN\}$ \ are independent copies of \ $V_{n,-1}$,
 \ where \ $(V_{k,0})_{k\geq-1}$ \ and \ $(V_{k,-1})_{k\geq-1}$ \ are second-order Galton--Watson
 processes (without immigration) with initial values \ $V_{0,0} = 1$,
\ $V_{-1,0} = 0$, \ $V_{0,-1} = 0$ \ and \ $V_{-1,-1} = 1$, \ and with the same offspring
 distributions as \ $(X_k)_{k\geq-1}$.

Moreover, if \ $(X_n)_{n\geq -1}$ \ is a second-order Galton--Watson process with immigration, then
 for each \ $n \in \NN$, \ we have
 \begin{equation}\label{2GWI_additive}
  X_n = V_0^{(n)}(X_0, X_{-1}) + \sum_{i=1}^n V_i^{(n-i)}(\vare_i, 0) ,
 \end{equation}
 where \ $\bigl\{V_0^{(n)}(X_0, X_{-1}), V_i^{(n-i)}(\vare_i, 0) : i \in \{1, \ldots, n\}\bigr\}$
 \ are independent random variables such that \ $V_0^{(n)}(X_0, X_{-1})$ \ represents the number of
 newborns at time \ $n$, \ resulting from the initial individuals \ $X_0$ \ at time \ $0$ \ and
 \ $X_{-1}$ \ at time \ $-1$, \ and for each \ $i \in \{1, \ldots, n\}$, \ $V_i^{(n-i)}(\vare_i, 0)$
 \ represents the number of newborns at time \ $n$, \ resulting from the immigration \ $\vare_i$ \ at
 time \ $i$.
\ Clearly, \ $(V_0^{(k)}(X_0, X_{-1}))_{k\geq-1}$ \ and \ $(V_i^{(k)}(\vare_i, 0))_{k\geq-1}$,
 \ $i \in \{1, \ldots, n\}$, \ are second-order Galton--Watson processes (without immigration) with
 initial values \ $V_0^{(0)}(X_0, X_{-1}) = X_0$, \ $V_0^{(-1)}(X_0, X_{-1}) = X_{-1}$,
 \ $V_i^{(0)}(\vare_i, 0) = \vare_i$ \ and \ $V_i^{(-1)}(\vare_i, 0) = 0$,
 \ $i \in \{1, \ldots, n\}$, \ and with the same offspring distributions as \ $(X_k)_{k\geq-1}$.

\section{Regularly varying functions}
\label{App1}

First, we recall the notion of regularly varying non-negative random variables.

\begin{Def}
A non-negative random variable \ $X$ \ is called regularly varying with index \ $\alpha \in \RR_+$
 \ if for all \ $q \in \RR_{++}$,
 \[
   \lim_{x\to\infty} \frac{\PP(X > qx)}{\PP(X > x)} = q^{-\alpha} .
 \]
\end{Def}

\begin{Lem}\label{rvexp}
If \ $X$ \ is a non-negative regularly varying random variable with index \ $\alpha \in \RR_{++}$,
 \ then \ $\EE(X^\beta) < \infty$ \ for all \ $\beta \in (-\infty, \alpha)$ \ and
 \ $\EE(X^\beta) = \infty$ \ for all \ $\beta \in (\alpha, \infty)$.
\end{Lem}
For Lemma \ref{rvexp}, see, e.g., Embrechts et al.\ \cite[Proposition A3.8]{EmbKluMik}.

\begin{Lem}\label{exposv}
 If \ $X$ \ and \ $Y$ \ are non-negative random variables such that \ $X$ \ is regularly varying with
 index \ $\alpha \in \RR_+$ \ and there exists \ $r \in (\alpha, \infty)$ \ with
 \ $\EE(Y^r) < \infty$, \ then \ $\PP(Y > x) = \oo(\PP(X > x))$ \ as \ $x \to \infty$.
\end{Lem}
For Lemma \ref{exposv}, see, e.g., Barczy et al.\ \cite[Lemma C.6]{BarBosPap}.

\begin{Lem}\label{svsv}
If \ $X_1$ \ and \ $X_2$ \ are non-negative regularly varying random variables with index
 \ $\alpha_1 \in \RR_+$ \ and \ $\alpha_2 \in \RR_+$, \ respectively, such that
 \ $\alpha_1 < \alpha_2$, \ then \ $\PP(X_2 > x) = \oo(\PP(X_1 > x))$ \ as \ $x \to \infty$.
\end{Lem}
For a proof of Lemma \ref{svsv}, see, e.g., Barczy et al.\ \cite[Lemma C.7]{BarBosPap}.

\begin{Lem}[Convolution property]\label{Lem_conv}
If \ $X_1$ \ and \ $X_2$ \ are non-negative random variables such that \ $X_1$ \ is regularly
varying with index \ $\alpha \in \RR_+$ \ and \ there exists \ $r \in (\alpha, \infty)$ \ with
\ $\EE(X_2^r) < \infty$, \ then \ $\PP(X_1 + X_2 > x) \sim \PP(X_1 > x)$ \ as \ $x \to \infty$,
\ and hence \ $X_1 + X_2$ \ is regularly varying with index \ $\alpha$.

If \ $X_1$ \ and \ $X_2$ \ are independent non-negative regularly varying random variables with
index \ $\alpha \in \RR_+$, \ then \ $\PP(X_1 + X_2 > x) \sim \PP(X_1 > x) + \PP(X_2 > x)$ \ as
\ $x \to \infty$, \ and hence \ $X_1 + X_2$ \ is regularly varying with index \ $\alpha$.
\end{Lem}
The statements of Lemma \ref{Lem_conv} follow, e.g., from parts 1 and 3 of Lemma B.6.1 of
 Buraczewski et al.\ \cite{BurDamMik} and Lemma \ref{svsv} together with the fact that the sum of two
 slowly varying functions is slowly varying.

\section{Regularly varying random sums}\label{rvsums}

Now, we recall conditions under which a random sum is regularly varying.
Combining Proposition 4.3 in Fa\"y et al.\ \cite{GilGAMikSam} and Theorem 3.2 in Robert and Segers
 \cite{RobSeg}, we obtain the following result (see, e.g., Barczy et al \cite[Lemma D.3]{BarBosPap}).
 
\begin{Pro}\label{FGAMSRS}
Let \ $\tau$ \ be a non-negative integer-valued random variable and let
 \ $\{\zeta, \zeta_i : i \in \NN\}$ \ be independent and identically distributed non-negative random
 variables, independent of \ $\tau$, \ such that \ $\tau$ \ is regularly varying with index
 \ $\beta \in \RR_+$ \ and \ $\EE(\zeta) \in \RR_{++}$.
\ In case of \ $\beta \in [1,\infty)$, \ assume additionally that there exists
 \ $r \in (\beta, \infty)$ \ with \ $\EE(\zeta^r) < \infty$.
\ Then we have
 \[
   \PP\biggl(\sum_{i=1}^\tau \zeta_i > x\biggr)
   \sim \PP\biggl(\tau > \frac{x}{\EE(\zeta)}\biggr)
   \sim (\EE(\zeta))^\beta \PP(\tau > x) \qquad \text{as \ $x \to \infty$,}
 \]
 and hence \ $\sum_{i=1}^\tau \zeta_i$ \ is also regularly varying with index \ $\beta$.
\end{Pro}

The next proposition is a consequence of part (ii) of Theorem 1 in Denisov et al.\ \cite{DenFosKor}
 (see, e.g., Barczy et al. \cite[Lemma D.5]{BarBosPap}).
 
\begin{Pro}\label{DFK_corr0}
Let \ $\tau$ \ be a non-negative integer-valued random variable and let
 \ $\{\zeta, \zeta_i : i \in \NN\}$ \ be independent and identically distributed non-negative random
 variables, independent of \ $\tau$, \ such that \ $\zeta$ \ is regularly varying with index
 \ $\alpha \in [1, \infty)$, \ $\PP(\tau = 0) < 1$ \ and there exists \ $r \in (\alpha, \infty)$
 \ with \ $\EE(\tau^r) < \infty$. 
\ In case of \ $\alpha = 1$, \ assume additionally that \ $\EE(\zeta) \in \RR_{++}$.
\ Then we have
 \[
   \PP\biggl(\sum_{i=1}^\tau \zeta_i > x\biggr) \sim \EE(\tau) \PP(\zeta > x) \qquad
   \text{as \ $x \to \infty$,}
 \]
 and hence \ $\sum_{i=1}^\tau \zeta_i$ \ is also regularly varying with index \ $\alpha$.
\end{Pro}

The next proposition is a consequence of Theorem 7 in Denisov et al.\ \cite{DenFosKor} (see, e.g.,
 Barczy et al. \cite[Lemma D.7]{BarBosPap}).

\begin{Pro}\label{DFK_corr}
Let \ $\tau$ \ be a non-negative integer-valued random variable and let
 \ $\{\zeta, \zeta_i : i \in \NN\}$ \ be independent and identically distributed non-negative random
 variables, independent of \ $\tau$, \ such that \ $\tau$ \ and \ $\zeta$ \ are regularly varying with
 index \ $\beta \in [1, \infty)$, \ and \ $\PP(\zeta > x) = \OO(\PP(\tau > x))$ \ as \ $x \to \infty$.
\ In case of \ $\beta = 1$, \ assume additionally that \ $\EE(\tau) \in \RR_{++}$ \ and
 \ $\EE(\zeta) \in \RR_{++}$.
\ Then we have
 \[
   \PP\biggl(\sum_{i=1}^\tau \zeta_i > x\biggr)
   \sim \EE(\tau) \PP(\zeta > x) + (\EE(\zeta))^\beta \PP(\tau > x) \qquad \text{as \ $x \to \infty$,}
 \]
 and hence \ $\sum_{i=1}^\tau \zeta_i$ \ is also regularly varying with index \ $\beta$.
\end{Pro}

\section{Moment estimations}
\label{App_moments}

We present an auxiliary lemma on the higher moments of second-order Galton--Watson processes.

\begin{Lem}\label{Lem_seged_momentr}
Let \ $(X_n)_{n\geq-1}$ \ be a second-order Galton--Watson process (without immigration) such that
 \ $\EE(X_0^r) < \infty$, \ $\EE(X_{-1}^r) < \infty$, \ $\EE(\xi^r) < \infty$ \ and
 \ $\EE(\eta^r) < \infty$ \ with some \ $r \in (1, \infty)$.
\ Then \ $\EE(X_n^r) < \infty$ \ for all \ $n \in \NN$.
\end{Lem}

\noindent{\bf Proof.}
We proceed by induction on \ $n$. 
\ For \ $n = 0$ \ and \ $n = -1$ \ the statement follows from the assumptions. 
Now let us assume that the statement holds for \ $1, \ldots, n - 1$.
\ For \ $n$, \ by power means inequality, we have
 \begin{align*}
  \EE(X_n^r \mid \cF_{n-1})
  &= \EE\left(\left(\sum_{i=1}^{X_{n-1}} \xi_{n,i}
                    + \sum_{j=1}^{X_{n-2}} \eta_{n,j}\right)^r \mid \cF_{n-1}\right) \\
  &\leq 2^{r-1}
        \EE\left(\left(\sum_{i=1}^{X_{n-1}} \xi_{n,i}\right)^r
                 + \left(\sum_{j=1}^{X_{n-2}} \eta_{n,j}\right)^r \mid \cF_{n-1}\right) \\
  &\leq 2^{r-1}
        \EE\left(X_{n-1}^{r-1} \sum_{i=1}^{X_{n-1}} \xi_{n,i}^r
                 + X_{n-2}^{r-1} \sum_{j=1}^{X_{n-2}} \eta_{n,j}^r \mid \cF_{n-1}\right) \\
  &= 2^{r-1}
     \left(X_{n-1}^r \EE(\xi^r) + X_{n-2}^r \EE(\eta^r)\right)
   < \infty
 \end{align*}
 for all \ $n \in \NN$, and hence
 \ $\EE(X_n^r) \leq 2^{r-1} \EE(X_{n-1}^r) \EE(\xi^r) +\EE( X_{n-2}^r) \EE(\eta^r) < \infty$
 \ by the induction hypothesis.
\proofend


\begin{thebibliography}{99}

\bibitem{KasDes}
Kashikar, A.S.; Deshmukh, S.R.
Probabilistic properties of second order branching process.
Ann.\ Inst.\ Statist.\ Math.
\textbf{2015}, \textsl{67(3)}, 557--572. 
DOI: 10.1007/s10463-014-0462-0.

\bibitem{KasDes2}
Kashikar, A.S.; Deshmukh, S.R.
Estimation in second order branching processes with application to swine flu data.
Comm.\ Statist.\ Theory Methods
\textbf{2016}, \textsl{45(4)}, 1031--1046.
DOI: 10.1080/03610926.2013.853796.

\bibitem{Kas}
Kashikar, A.S.
Estimation of growth rate in second order branching process.
J. Statist.\ Plann.\ Inference
\textbf{2017}, \textsl{191}, 1--12.
DOI: 10.1016/j.jspi.2017.06.003.

\bibitem{Lat}
Latour, A.
The multivariate GINAR(p) process.
Adv.\ in Appl.\ Probab.
\textbf{1997}, \textsl{29(1)}, 228--248.
DOI: 10.2307/1427868.

\bibitem{DuLi}
Du, J.-G.; Li, Y.
The integer valued autoregressive (INAR($p$)) model.
J. Time Series Anal. 
\textbf{1991}, \textsl{12(2)}, 129--142.
DOI: 10.1111/j.1467-9892.1991.tb00073.x.

\bibitem{BarBosPap}
Barczy, M.; B\H osze, Zs.; Pap, G.
Regularly varying non-stationary Galton--Watson processes with immigration.
\textbf{2018+}, http://arxiv.org/abs/1801.04002.

\bibitem{BarBosPap2}
Barczy, M., B\H osze, Zs.; Pap, G.
On tail behaviour of stationary second-order Galton--Watson processes with immigration.
\textbf{2018+}, https://arxiv.org/abs/1801.07931.

\bibitem{AthNey}
Athreya, K.B.; Ney, P.E.
\textsl{Branching Processes};
Springer-Verlag: New York-Heidelberg, 1972.

\bibitem{EmbKluMik}
Embrechts, P., Kl\"uppelberg, C.; Mikosch, T.
\textsl{Modelling Extremal Events for Insurance and Finance};
Springer: Berlin, 1997.

\bibitem{BurDamMik}
Buraczewski, D.; Damek, E.; Mikosch, T.
\textsl{Stochastic models with power-law tails. The equation $X = A X + B$};
Springer: Cham, 2016.

\bibitem{GilGAMikSam}
Fa\"y, G.; Gonz\'alez-Ar\'evalo, B.; Mikosch, T.; Samorodnitsky, G.
Modeling teletraffic arrivals by a Poisson cluster process.
Queueing Syst.
\textbf{2006}, \textsl{54}, 121--140.
DOI: 10.1007/s11134-006-9348-z.

\bibitem{RobSeg}
Robert, C.Y.; Segers, J.
Tails of random sums of a heavy-tailed number of light-tailed terms.
Insurance Math.\ Econom.
\textbf{2008}, \textsl{43}, 85--92. 
DOI: 10.1016/j.insmatheco.2007.10.001.

\bibitem{DenFosKor}
Denisov, D.; Foss, S.; Korshunov, D.
Asymptotics of randomly stopped sums in the presence of heavy tails.
Bernoulli
\textbf{2010}, \textsl{16(4)}, 971--994.
DOI: 10.3150/10-BEJ251.

\end{thebibliography}
\end{document}